\documentclass{amsart}
\usepackage{graphicx}
\usepackage[colorlinks, linkcolor=blue, anchorcolor=blue, citecolor=blue]{hyperref}
\usepackage{amsfonts,amssymb,latexsym,epsfig,amsmath,enumerate}
\usepackage{filecontents}
\newtheorem{thm}{Theorem}[section]

\newtheorem{lem}[thm]{Lemma}
\newtheorem{prop}[thm]{Proposition}
\theoremstyle{definition}
\newtheorem{defn}[thm]{Definition}
\theoremstyle{remark}
\newtheorem{rem}[thm]{Remark}
\newtheorem{exam}[thm]{Example}

\numberwithin{equation}{section}
\newcommand{\norm}[1]{\left\Vert#1\right\Vert}
\newcommand{\abs}[1]{\left\vert#1\right\vert}
\newcommand{\set}[1]{\left\{#1\right\}}
\newcommand{\Real}{\mathbb R}
\newcommand{\Epc}{\mathbb E\,}
\newcommand{\eps}{\varepsilon}

\newcommand{\R}{\mathbb{R}}
\newcommand{\N}{\mathbb{N}}
\newcommand{\Rd}{\mathbb{R}^d}
\newcommand{\supp}{\text{supp }}

\DeclareMathOperator{\essinf}{ess\,inf\,}

\begin{document}

\title[On solving SDEs by density approximations]
{On numerical density approximations of solutions of SDEs with
  unbounded coefficients}


\author[L. Chen]{Linghua Chen}%
\address{Norwegian University of Science and Technology, NO-7491, Trondheim, Norway}%
\email{linghuac@math.ntnu.no}%

\author[E. R. Jakobsen]{Espen Robstad Jakobsen}%
\address{Norwegian University of Science and Technology, NO-7491, Trondheim, Norway}%
\email{erj@math.ntnu.no}%

\author[A. Naess]{Arvid Naess}%
\address{Norwegian University of Science and Technology, NO-7491, Trondheim, Norway}%
\email{arvidn@math.ntnu.no}%

\thanks{E. R. Jakobsen is partially supported by the
NFR Toppforsk project 250070, Waves and Nonlinear Phenomena (WaNP).
The authors are thankful to Harald Hanche-Olsen and
Peter Lindqvist for their valuable inputs to Proposition \ref{prop dissipative}
and Lemma \ref{lem ell reg}, respectively.}%

\subjclass{60H35, 65M12, 47D07}%
\keywords{Stochastic Differential Equations, Numerical Method,
    Path Integration, Density Tracking, Probability Density,  Semigroup
    Generation, Convergence}%

\date{\today}%

\begin{abstract}
  We study a numerical method to compute probability
  density functions of solutions of stochastic differential
  equations. The method is sometimes called the numerical path integration
  method and has been shown to be fast and accurate in application
  oriented fields. In this paper we provide a rigorous analysis of
  the method that covers systems of equations with unbounded coefficients.
Working in a natural space for densities, $L^1$, we obtain
stability, consistency, and new convergence results for the method, new
well-posedness and
  semigroup generation results for the related Fokker-Planck-Kolmogorov
  equation, and a new and rigorous connection to the corresponding probability
  density functions for both the approximate and the exact
  problems.
To prove the results we combine semigroup and PDE arguments in a new
way that should be of independent interest.
\end{abstract}
\maketitle

\section{Introduction}
Over the past decades there has been a large number of publications in
the field of stochastic dynamics and its various application
areas --
including physics, biology, engineering, and finance
\cite{MK2004,applebaum2004levy,tankov2003financial}.
In this field, the response
of dynamical systems to stochastic excitation is studied, and the
typical model
is a (system of) stochastic differential equations
(SDEs) of the form
\[
  \begin{cases}
      dY_t = b(t, Y_t)\, dt + \sigma(t, Y_t) \, dB_t, \\
      Y_0 = Z,
  \end{cases}
\]
where $b: \Real_+\times \Real^d \to \Real^d$, $\sigma: \Real_+\times
\Real^d \to \Real^{d\times n}$, $B_t$ is an $n$-dimensional Brownian
motion, and the initial data $Z$ is a random variable in $\R^d$
independent of $B_t$. The solution $Y_t$ of the SDE is a
 state space process in $\R^d$.

In most cases the solutions of such problems must be computed
numerically, and various discrete  approximation methods are widely
used in many application areas \cite{kloeden2013numerical}.
There are two main approaches: (i) Path-wise approximations of the SDE
based on stochastic simulation,
and (ii) approximations of the statistics or distributions of the
SDE. The first approach is more efficient in high dimensions and the
second in low dimensions. For path-wise approximations we refer to
\cite{kloeden2013numerical} for the classical literature and e.g. to
\cite{LiLy2012} and references
therein for some promising recent developments. Approach (ii) is
deterministic and based on approximating
the forward or backward Kolmogorov partial differential equations.
In this paper we study a method of the second type, called
the numerical path integration method \cite{Naess199391, linetsky1998path, Rosa-clot99apath,Yu20041493,
  NaessSkaug2007, Chai2015} or density tracking method
\cite{BM16}. This method is an explicit deterministic iteration
scheme that produces approximate probability density functions (PDFs) for the
solution of the SDE. The iteration step is based on a short time
approximation of the SDE. For simplicity we use here the strong
Euler-Maruyama method, the most basic numerical method for SDEs. The
convergence of our path integration method is therefore equivalent to
the convergence of the PDFs of the Euler-Maruyama method.

The path integration/density tracking approach
(i.e. simulating the PDFs)  enjoys several favorable properties.
First, it introduces an extra perspective to the system,
which enables deeper insights and invites broader mathematical tools.
Secondly, as an explicit method, one can formally implement
the path integration algorithm on a vast number of scenarios.
Since the formulation is deterministic, it is also
free from perturbation by extreme outcomes during stochastic simulation.
Finally, the result of the method is an explicit density function
rather than bundle of random paths. This means that many
characteristics of the system become more transparent,
and can be captured and displayed by e.g. visualisation methods.

The numerical path integration method has been applied in many fields,
including financial mathematics
\cite{linetsky1998path, Rosa-clot99apath, NaessSkaug2007, Chai2015}.
Many of these studies show that it can provide highly
accurate numerical solutions \cite{Naess199391, Yu20041493, Mo2008phd,BM16}.
Even though convergence problems have been reported in
  some cases,
  cf. Section 7.1 in \cite{Mo2008phd}, little or no emphasis
has been devoted to conditions for convergence of the method in the path
integration literature.
A very natural and relevant mode of
  convergence for this method is strong
$L^1$-convergence of the resulting densities.
Such convergence seems not to be a direct consequence of either strong
or weak convergence of the Euler-Maruyama
scheme. In fact, $L^1$-convergence of densities implies weak
convergence of processes but the 
converse is not true in general. Nevertheless, very recently there are
some results for one-dimensional problems in \cite{BM16}, and when the
coefficients are bounded, $L^1$-convergence results follow from some of very
precise error expansions for densities of the Euler-Maruyama method,
e.g. in \cite{BallyTalayII1995, Guyon2006, GobetLa2008}, or from
so-called Feynman formulas in mathematical physics
\cite{Butko2016}. More on this below. However from an application
point of view, it is important to consider also models with unbounded
coefficients, since they naturally appear in many of the papers and
applications mentioned above.\footnote{By a change of variables, it is
often possible to reduce to a problem with bounded coefficients and
then use existing $L^1$ convergence results. However going
back to the original variables, we would then get a weaker form of
convergence, one that no longer implies $L^1$ convergence.}

The aim of this paper is to provide a 
  rigorous analysis of method in the path integration setting
  along with the $L^1$-convergence of the PDFs it produces in the case
  of unbounded coefficients.
Specifically, we obtain stability, consistency,
  and new convergence results for the method, new well-posedness
  and semigroup generation results
 for the associated Fokker-Planck-Kolmogorov (FPK)
  equation, and a new and rigorous connection to the corresponding probability
  density functions for both the approximate and exact
  problems.
In particular, our $L^1$ semigroup generation result for our rather
general elliptic FKP operator 
seems to be new.

We also investigate in details the method applied to the
Ornstein-Uhlenbeck process which has unbounded drift. Using explicit
transition densities and direct computations, we show that the method
(i) converge in the general case and (ii) converge with a linear rate
if additional moment conditions are assumed. The second result is
consistent with results for bounded coefficients in e.g.
\cite{BallyTalayII1995, Guyon2006, GobetLa2008}.
Even though we do not have a counter example, the very
        explicit proof seems to indicate that convergence rates can
        only be obtained under some additional moment assumptions on the densities.

We study the numerical path integration method from a PDE and semigroup point
of view, an approach which seems to be new to this setting at least
when it comes to the convergence results. The main idea is to show
that (i) the FPK operator corresponding to the SDE generates a
semigroup, (ii) the iteration of
the discrete path integration operator does converge to this
semigroup  (the solution of the
FPK equation), and finally (iii) that convergence of the densities follows
from these results. For (i), we prove that the FKP operator and its
adjoint are dissipative and then use the Lumer-Phillips theorem to prove
that the FPK operator generates a
strongly continuous contraction semigroup on $L^1(\Real^d)$.
To prove the convergence results in (ii), we use the Chernoff product
formula which is a generalisation of the well-known Kato-Trotter
product formula. Note that our approach to (ii) is similar to the
approach for Feynman formulas in \cite{Butko2016}.

A main difficulty is then to prove the semigroup generation in our
setting of unbounded coefficients. Classical results usually assume
the coefficients to be bounded, cf. the discussion in Chapter 11 of
\cite{lasota2013chaos}.
Further complications comes from the fact that we have to  work in
the space $L^1$ which is
not reflexive and the adjoint space $L^{\infty}$ which is not separable.
In fact one of the most difficult parts of our paper is
to show that the adjoint operator is dissipative on $L^{\infty}$. To
do this we develop a new and non-trivial argument using ideas from viscosity
solution theory \cite{CIL}, the weak
Bony maximum principle \cite{Bony1967, Lions1983}, and recent elliptic
regularity results from \cite{ZhangBao2013}. This argument could be
of independent interest.

Other authors have obtained semigroup generation results for similar
problems with different techniques.
Yet, most work seems to be devoted to the generation problems for
the {\em adjoint} of the FPK operator in various situations.  For instance
\cite{DaPrato2001333} is concerned with the space of
bounded (continuous) functions and unbounded coefficients, while in
\cite{Angiuli2010} the authors give results in $L^1$ for bounded
coefficients.  The author in \cite{Cerrai2000151}
considered degenerate operators where the coefficients have
bounded first and second derivatives.
The result closest to our case seems to be found in
\cite{Fornaro2007747} where the drift term had
to be dominated by the "square root" of the diffusion term and no zero
order term is allowed. A semigroup generation result of the adjoint
operators on $L^1$ is proved there.
However, we cannot use the results in \cite{Fornaro2007747}
in this paper since our non-divergence form FPK operators are
equivalent to divergence form operators with additional unbounded zero
order terms.

In addition to providing rigorous results for the
numerical path integration method, this paper seems to represent the
first attempt to use semigroup and PDE arguments to obtain strong
convergence results for probability densities of SDE approximations
(but see also \cite{Butko2016}).
It is an interesting question how far such methods can be pushed compared
to other methods.
As far as we know, two other methods are
described in the literature: A parametrix approach of Konakov et
al. \cite{KonakovM2000, konakov2010, MenozziLe2010},
 and an approach based on
Malliavin calculus developed by Talay, Bally, et al., see e.g.
\cite{BallyTalayII1995, Guyon2006, Bally2014} and references therein.
In the parametrix approach, a
weighted $L^\infty$-estimate on the convergence of (transition
probability) densities,
as well as non-asymptotic bounds for
the densities \cite{MenozziLe2010},
is obtained under uniform ellipticity
assumption and boundedness of coefficients.
In the Malliavin approach various error expansions
and estimates are obtained under boundedness, certain smoothness, and
ellipticity or H{\"o}rmander conditions.

Under additional assumptions, the $L^1$-convergence of
densities of the Euler-Maruyama scheme can be obtained from
some of these Malliavin results.
In \cite{BallyTalayII1995} a linear convergence rate was obtained for
SDEs with smooth and bounded coefficients\footnote{The coefficients are bounded in
  \cite{BallyTalayII1995}. This follows from their condition (H), estimate
  (6), and footnote 1. In particular, in \cite{BallyTalayII1995} the
  number $\gamma_0$, introduced by (1.10) in \cite{KusuokaStr1985}, is zero.},
but the behavior of the estimate is unclear for small time.
In \cite{Guyon2006} the estimate was improved but
it would still blow up as $t\to 0$.
Later the authors of \cite{GobetLa2008}
derived very accurate estimates that would lead to a linear rate
which is uniform for small time. At the same time they
lowered the regularity requirements of the coefficients to $C^{1,3}_b$.
For the Euler-Maruyama method, the estimates of
\cite{BallyTalayII1995, Guyon2006, GobetLa2008} requires bounded
coefficients and  uniform
ellipticity. According to Remark 2.22 of \cite{Bally2014}, the
ellipticity requirement can be relaxed to a so-called weak
H{\"o}rmander condition when coefficients belong to $C_b^\infty$. In this case,
the results of \cite{Bally2014} give $L^1$-convergence for fixed
times (but no estimates) of densities/PDFs for the Euler-Maruyama scheme.

To summarise, in our approach to study properties and convergence of the
numerical path integration method we obtain the following results: New
semigroup
generation and well-posedness results for
the FPK equation, allowing for coefficients with linear growth and
$C^4$ differentiability, no invariant measure is required,
the growth of drift and diffusion terms are not necessarily related to each
another, and the initial data is not required to have any finite moments
or differentiability. Under the same conditions,
we prove $L^1$ convergence, locally uniformly in time (so no deterioration for small time).
To our knowledge, this seems to be the first general strong $L^1$-convergence
result for densities/PDFs of SDEs with unbounded coefficients.
Even though we do not obtain any rate of convergence or error
expansion in general, by studying the Ornstein-Uhlenbeck process,
we find indications that no error estimate may exist without
additional (moment) assumptions on the densities/PDFs. Hence when the
coefficients are unbounded, the best result for general densities
could be mere convergence without any rate.

We also emphasise that our results connects the mild solution of the
 FPK equation and the
PDF of the SDE, and that existence and convergence of the latter
 follows from the existence and  convergence (by semigroup
methods) of the
former. More refined existence
results for PDFs have been obtained by probabilistic methods, we refer
to e.g. \cite{Bouleau1986, fournier2010, Hayashi2012, HayashiKHY2014}
and references therein.

\subsection*{Layout and Notation}
The rest of this paper is laid out as follows.
In Section \ref{sec main results} we state our assumptions and the main
results. Included is a discussion of the connection between SDEs, FPK
equations and densities and the definition of the path integration method.
At the end of the section, we discuss our results and give some examples.
The rest of the paper is then devoted to the proofs of these results.
In Section \ref{sec pf convg of disc pi} we prove
the well-posedness of the FPK equation, the connection to probability
densities, and the convergence result -- assuming the semi-group
generation and a strong $L^1$-consistency result for the path
integration method,
which are proved in Section \ref{Pf2} and
 Section \ref{sec gen result}, respectively.
At the end, Section \ref{sec rate} is an appendix to the calculation
of the convergence rate of Example \ref{exam_lgv}.

The following notation and abbreviations are used
$\norm{\cdot}_1:= \norm{\cdot}_{L^1(\Real^d)}$,
$\norm{\cdot}_{\infty}:= \norm{\cdot}_{L^{\infty}(\Real^d)}$,
$\partial_t:= \frac{\partial}{\partial t}$,
$\nabla:=  \left(\frac{\partial}{\partial x_1},\cdots,\frac{\partial}{\partial x_d}\right)^T
            =:(\partial_1 ,\cdots,\partial_d )^T$,
$\essinf$ is the
essential infimum, $\Epc$ denotes the expectation; $C_b^k$,
$C_c^\infty$, $\mathcal D'$ the spaces of functions with bounded continuous
derivatives up to $k-$th order, smooth compactly supported
functions, and distributions ($\mathcal D'$ is the dual of
$C_c^\infty$) respectively,
$C\big([0,T];L^1(\R^d)\big)$ the space of functions $u(x,t)$ such
that  $\sup_{t\in[0,T]}\|u(\cdot,t)\|_1<\infty$,
and for all $t\in[0,T]$, $\|u(\cdot,s)-u(\cdot,t)\|_1\to 0$ as $s\to t$.
PDF -- probability density
function, SDE -- stochastic differential equation, FPK -- Fokker-Planck-Kolmogorov.

\section{Main Results}\label{sec main results}
In this section we formulate the discrete path integration method, state
the assumptions and main results, and provide a discussion and examples.
We restrict ourselves to time-homogeneous SDEs,
i.e.  equations of the form
\begin{equation}\label{eq_general_sde}
    \begin{cases}
        d Y_t = b(Y_t)\, dt + \sigma(Y_t) \, dB_t,\\
        Y_0 = Z,
    \end{cases}
\end{equation}
where $b = (b_1,\cdots,b_d)^T:\R^d\to\R^d$ and
$\sigma=(\sigma_{ij})_{d\times n}:\R^d\to\R^{d\times n}$ are
functions, $B_t$ a $d$-dimensional Brownian motion,
and $Z$ is a random variable in $\R^d$.

Under suitable assumptions (cf. \cite{oksendal2010stochastic} or
\cite{protter2005stochastic}), the solution
$Y_t$ of \eqref{eq_general_sde},
is a Markov process, with infinitesimal generator $A^*$ defined as
\begin{equation}\label{eq adj op on Cb}
    A^*\phi:=\sum_{i=1}^d
    b_i(x)\partial_i\phi(x)+\frac{1}{2}\sum_{i,j=1}^d
    a_{ij}(x) \partial_i\partial_j\phi(x)
    \qquad\text{for all}\qquad \phi\in C_c^{\infty}(\Real^d),
\end{equation}
where
$$a:=\sigma\sigma^T=:(a_{ij})_{d\times d}.$$
Moreover,
$Y_t$ has a PDF $u(t,x)$ -- a non-negative $L^1$ function $u$ such that
\begin{equation}\label{eq def of PDF}
    \Epc \phi(Y_t)=\int_{\Real^d} \phi(x) u(t,x)dx\quad
    \text{ for all }\phi\in C_b^0(\Real^d).
\end{equation}
Via an adjoint argument, $u$ formally satisfies the FPK equation
\begin{equation}\label{eq FPK}
    \begin{cases}
        \partial_t u(t,x)  =Au(t,x),\qquad& x\in\R^d,\ t>0,\\[0.2cm]
        u(0,x)  = u_0(x),&x\in\R^d,
    \end{cases}
\end{equation}
where $u_0$ is the PDF of $Z$ and $A$ the adjoint of $A^*$,
\begin{equation}\label{eq def of A}
    A v(x):= -\sum_{i=1}^d\partial_i\big(b_i(x)v(x)\big)
        +\frac{1}{2}\sum_{i,j=1}^d\partial_i\partial_j\big(a_{ij}(x)v(x)\big).
\end{equation}
The evolution of the $u$ can be regarded as the action of the
propagation operator (or semigroup) $P_s$ defined by
\begin{equation}\label{eq_def_PI}
    u(y,t+s)=[P_{s}u(\cdot \,,t)](y)\qquad\text{for all}\qquad
     t,s>0,\ y\in\R^d.
\end{equation}

Let us state our assumptions.
\medskip
\begin{itemize}
  \item[\bf (C1)]
    $b:\R^d\to\R^d$ and $\sigma:\R^d\to\R^{d\times n}$ are $C^4$
    functions, and there exists a constant $K>0$ such that for all
    $x\in \Real^d$, $j=1,\dots,n$, and $i,k=1,\dots,d$,
    \[
        |\partial_k\sigma_{ij}(x)|+|\partial_k b_i(x)|  \leq K.
    \]
  \item[\bf(C2)]
   There exists $\alpha>0$ such that for all $x, y\in \Real^d$,
    \[
        y^T a(x) y \geq \alpha |y|^2.
    \]
  \item[\bf(C3)]
    $Z$ and $B_t$ are independent,
    and $Z$ has a PDF $u_0$, i.e. $u_0\geq 0$ and $\norm{u_0}_1=1$.
\end{itemize}
\begin{rem}
Standard strong well-posedness of the SDE  (see below) requires only Lipschitz
coefficients. The additional regularity and the uniform
ellipticity condition (C2) is needed to study the PDFs of the SDEs and
path integration method (cf. proofs of Proposition
    \ref{prop lim of disc PI} and Lemma \ref{lem critical seq of bounded function}).
Similar but stronger assumptions, including bounded coefficients, are
used by all other papers discussing densities of SDE approximations, see
 the introduction and e.g. \cite{GobetLa2008}.
\end{rem}

The following result is then classical, cf. Theorem V.7 in \cite{protter2005stochastic}.
\begin{prop}\label{prop SDE exist and uniq}
    Assume (C1) and (C3). Then the SDE \eqref{eq_general_sde} has
    a unique strong solution $Y_t$.
\end{prop}

Now we give our first main result -- well-posedness for the FPK
equation \eqref{eq FPK}.
\begin{thm}\label{thm solve FPK}
    Assume (C1), (C2), and $u_0 \in L^1(\Real^d)$. Then the FPK
    equation \eqref{eq FPK}
    has a unique (mild) solution $u\in C([0,T];L^1(\Real^d))$.
\end{thm}

The derivation of FPK equation \eqref{eq FPK} was formal.
Yet our second main result confirms that the probability distribution of
\eqref{eq_general_sde} coincides with the (mild) solution of FPK
equation \eqref{eq FPK}.
\begin{thm}\label{thm FPK PDF}
    Assume (C1) -- (C3), and let $Y_t$ and $u(t,x)$
    be the solutions of \eqref{eq_general_sde}
    and \eqref{eq FPK} respectively.
    Then \eqref{eq def of PDF} holds and $u(t,\cdot)$ is the PDF of $Y_t$.
\end{thm}

The proofs of Theorems \ref{thm solve FPK} and \ref{thm FPK PDF} will
be given in the next section.

\begin{rem}
    \noindent $(a)$\quad The proof shows that Theorem \ref{thm solve
      FPK} still holds if the $C^4$ regularity in assumption (C1)
    is reduced to $b\in C^1$ and $\sigma\in C^2$.
    \smallskip

    \noindent $(b)$\quad A difference between Theorem \ref{thm solve FPK}
    and existing results is that we allow (and need to allow) the zero order
    coefficient to be unbounded, cf. e.g. Section 4 in \cite{Fornaro2007747}.
    \smallskip

    \noindent $(c)$\quad The solution $u$ we obtain in
    Theorem \ref{thm solve FPK} is a mild solution of \eqref{eq FPK}
    (cf. \cite{engel2006one}). In fact we have more regularity than
    stated since $u(\cdot,t)$ belongs to the domain of $A$. A limit
    argument immediately shows that $u$ is also a distributional
    solution of \eqref{eq FPK}. In view of uniform
    ellipticity of $A$ and smoothness of its coefficients, it is standard
    to prove higher regularity and that $u$ is a classical solution.
\end{rem}

We now introduce and analyse the discrete path integration method. It
is a method for
computing approximations of the PDF of solution process of
\eqref{eq_general_sde} by:
\begin{enumerate}[($i$)]
    \item Finding a short-time approximation $\bar P_\tau$ of the propagator/semigroup $P_\tau$, and
        \smallskip
    \item computing the approximation $\bar u_n$ at time
        $t=n\tau$ using the explicit iteration
        \begin{equation}\label{PI}
            \bar u_{n+1}(x)=(\bar P_\tau\bar u_{n})(x)\quad \text{for}\quad n=0,1,\dots.
        \end{equation}
\end{enumerate}

The most obvious way to find a $\bar P_\tau$ is through a
time-discretization of the underlying SDE. For simplicity, we consider
here the Euler-Maruyama method \cite{kloeden2013numerical}, the most basic
and widely used SDE approximation scheme: Fix a time step $\tau>0$ and
let $\Delta B_{n}:=B_{n \tau} - B_{(n-1)\tau}$, and define the Markov
chain  approximation $X_n$ of the solution process $Y_{t}$  of the
SDE \eqref{eq_general_sde} by
\begin{equation}\label{eq EM disc}
    \begin{cases}
        X_{0} = Z,\\[0.2cm]
        X_{n+1} = X_n + b(X_n) \tau +
            \sigma(X_n)\Delta B_{n},\qquad n=0,1,\dots.
    \end{cases}
\end{equation}

In the rest of the paper we will use the following notation.
\begin{defn}\label{def1}\
    \noindent $(i)$\quad $\bar{u}_n$ is the PDF of the solution $X_n$ of
      \eqref{eq EM disc}.
    \smallskip

    \noindent $(ii)$\quad $\bar{P}_{\tau}$ is the one step propagator of
    $\bar{u}_n$, i.e.
    \begin{align*}
      \bar{u}_{n+1} (y) & = (\bar{P}_{\tau} \bar
        u_{n})(y)\qquad\text{for all}\qquad y\in\R^d,\
        \tau>0, \ n\geq 0.
    \end{align*}

    \noindent $(iii)$\quad $\bar k(y,x,\tau)$ is the transition kernel of
    $\bar{P}_{\tau}$ -- a non-negative $L^1$ function satisfying
    \begin{align}\label{Ptau}
         (\bar{P}_{\tau} \bar{u}_n)(y)=\int_{\R^d} \bar
                k(y,x,\tau)\bar{u}_n(x)dx\qquad\text{for all}\qquad y\in\R^d,\
                \tau>0, \ n\geq 0.
    \end{align}
\end{defn}

Since $\Delta B_{n}$ is a multivariate Gaussian variable, we have
the following classical result.
\begin{lem}\label{lem disc ker in F sp}
    Assume $a:=\sigma \sigma^T$ is strictly positive definite, and let
     $\bar k$ be given by Definition \ref{def1}. Then for all
     $x,y\in\R^d$ and $\tau>0$,
    \begin{gather}
    \mathcal{F}_y \bar{k}(\xi,x,\tau)  =  \exp\set{i\xi^T(x + b(x)
        \tau)-\frac{1}{2}\xi^T a(x)\xi},\nonumber \\
        \intertext{where $\mathcal{F}$ denote the Fourier transform, and}
        \bar{k}(y,x,\tau)=
                    \frac{\exp\set{-\frac{1}{2}(y-(x+b(x)\tau))^T (\tau a)^{-1}(x)(y-(x+b(x)\tau))}}
                    {(2\pi)^{2/d} |\det (\tau a(x))|^{1/2}}.\label{eq disc PI kernel}
    \end{gather}
\end{lem}

Existence and uniqueness of the solution $X_n$ come for free,
since the iteration \eqref{eq EM disc} is explicit. In view
of Lemma \ref{lem disc ker in F sp}, we then have the following result.

\begin{prop}[Well-posedness] Assume (C1) -- (C3).
    Then there exists a unique PDF $\bar u_n$ of the solution
    $X_n$ of \eqref{eq EM disc}. Moreover, $\bar u_n$ is explicitly given
    by \eqref{PI}, \eqref{Ptau}, and \eqref{eq disc PI kernel}.
\end{prop}

Note that (C2) implies that $a(x)$ is strictly positive definite.
A classical computation using \eqref{eq disc PI kernel} then shows
that $\|\bar k(\cdot,x,\tau)\|_1=1$, and hence (see the next Proposition) that
$\|\bar{P}_\tau\|= 1$. Proofs can be found e.g. in the discussion of Markov
operators in \cite{lasota2013chaos}.
\begin{prop}[$L^1-$stability]\label{prop disc PI contract}
    Assume (C2), $\tau>0$, $u\in
    L^1(\Real^d)$, and $\bar P_\tau$ is given by Definition
    \ref{def1}. Then
    $\|\bar{P}_{\tau}u\|_1\leq \|u\|_1$.
\end{prop}
By linearity, we immediately get continuous dependence on the initial data,
$$\|\bar{P}_{\tau}(u-v)\|_1\leq \|u-v\|_1,$$
and by iteration that the discrete path integration method
is $L^1$-{\em stable}:
$$\|\bar u_n\|_1\leq \|u_0\|_1\qquad\text{for all}\qquad n=1,2,\dots. $$

Next we study the $L^1$-{\em consistency} of the discrete scheme.
\begin{prop}[$L^1-$consistency]\label{prop lim of disc PI}
    Assume (C1), (C2), $\tau>0$, and $u\in C_c^{\infty}(\Real^d)$.
    Then there exists a constant $C>0$ such that
        \[
          \norm{\frac{\bar{P}_{\tau}u-u}{\tau}-Au}_1\leq C \tau.
        \]
\end{prop}
We prove this result in Section \ref{Pf2}. Higher regularity of $b$ or
$\sigma$ will not  improve the (linear) rate in Proposition \ref{prop lim of disc PI}.
It is the
maximal rate for the problem, cf. Proposition \ref{prop lgv rate} and its proof.

Finally we state the third main result of this paper, the
$L^1$-convergence of the discrete path integration method.
\begin{thm}[$L^1$-convergence]\label{thm disc PI converges}
    Assume (C1), (C2),  $u_0\in L^1(\Real^d)$, $u(t,x)$ is the
    solution of FPK equation \eqref{eq FPK}, and define
    $\bar{u}_n(t,x):=(\bar{P}_{t/n}^n u_0)(x)$. Then
    \begin{equation}\label{eq sg converg}
        \lim_{n\to \infty}\sup_{t\in [0,T]}
        \norm{u(\cdot,t)-\bar{u}_n(\cdot,t)}_1=0.
    \end{equation}
\end{thm}
\begin{rem}
    $(a)$\ Since the PDFs in general are  (non-negative)
    $L^1$-functions, strong $L^1$-convergence is a very natural mode of
    convergence to consider.
    \smallskip

    \noindent $(b)$\  By assumption the initial PDF $u_0$ is only required to belong to
    $L^1$. It is not required to have any finite moments or differentiability
    to obtain an $L^1$-convergence, and the convergence is uniform
    in time on $[0,T]$ for any $T>0$.
    \smallskip

    \noindent $(c)$\ $L^1$-convergence of PDFs is not a direct consequence
    of either strong or weak convergence of the solution process.
    It is strictly stronger
    than weak convergence of the corresponding
    process by Proposition \ref{L1weak} and Example \ref{exweak} in
    Appendix \ref{app:weak}.
\end{rem}

When the coefficients are also bounded, many authors have obtained not only
convergence but even error estimates and error expansions for
densities of the Euler-Maruyama scheme (see the introduction).

In the unbounded
coefficients case however, {\em there may not be any general $L^1$ error
bound without extra moment assumptions on the density}, and hence
 mere convergence would be the best one can hope for.
A first indication of
this appears already in the
Gaussian bounds on the derivatives of the SDE's transition probability
density $p(t,x,y)$ in \cite{KusuokaStr1985}.
When the coefficients are unbounded, these bounds
are no longer bounded in $x$, and hence can no
longer be used to derive the same estimates of approximations of $p(t,x,y)$ as
in the bounded case \cite{BallyTalayII1995,GobetLa2008}.
Another indication is given by Example \ref{exam_lgv} below. In this
example we study the Ornstein-Uhlenbeck process. We
find both the  exact and approximate transition kernels
and compute explicitly the conditions for having either mere
convergence in $L^1$ or a linear rate for the corresponding PDFs. In
the latter case, we find that additional moment assumptions are
needed (cf. Proposition \ref{prop lgv rate}).

\begin{exam}\label{exam_lgv}
    The Ornstein-Uhlenbeck process $Y_t$ is given by
    \begin{equation}\label{eq_lgv}
        d Y_t = b Y_t d t + \sigma d B_t  \qquad\text{and}\qquad Y_0=x,
    \end{equation}
    where $b,\;\sigma \in \Real$ and $Y_t,\,B_t$ are processes in
    $\R$. The unique strong solution is
    \[
        Y_t = e^{b t}\left(x + \sigma\int_0^t e^{-b s} d B_s \right).
    \]
    Due to the non-stochastic integrand in the It\^{o}-integral it is easy to see that
    this process has a Gaussian law (cf. e.g. Section 11.5
    in \cite{engel2006one}),
    \[
        N\left(xe^{bt},\frac{\sigma^2}{2b}\left(e^{2bt}-1\right)\right),
    \]
    and hence the transition kernel is given by
    \begin{equation}\label{conti_ker_lgv}
        k(y,x,t) = \frac{1}{\sqrt{2 \pi\frac{\sigma^2}{2b}\left( e^{2bt}-1\right)}}
        \exp\set{- \frac{\left(y-xe^{bt}\right)^2}{2\frac{\sigma^2}{2b}\left(e^{2bt}-1\right)}}.
    \end{equation}

    The corresponding Euler-Maruyama scheme with $\tau=t/n$ is
    \begin{align}\label{ex_scheme}
      \begin{cases}
          X_0 = x, \\
          X_{n+1} = X_n+b X_n \tau+\sigma\Delta B_{n},
      \end{cases}
    \end{align}
  and  has the following one step transition kernel by Lemma \ref{lem disc ker in F sp},
    \begin{equation}\label{disc_ker_lgv}
      \bar{k}(y,x,\tau) = \frac{1}{\sqrt{2 \pi \tau \sigma^2}}
        \exp\set{-\frac{(y - (1 + b \tau)x)^2}{2\tau \sigma^2}}.
    \end{equation}
    For $t=\tau$ small, it is clear that
    $\bar{k}(y,x,\tau)$ is close to $k(y,x,t)$.

For the convergence, we have the following result.
    \begin{prop}\label{prop lgv rate}
        Assume $b$, $\sigma\in \R$, $T>0$, and $u_0\in L^1(\R)$. Then
        $$\lim_{n\to\infty}\sup_{t\in [0,T]}\norm{P_tu_0-\bar P_{t/n}^n u_0}_1=0.$$
        If in addition, either
$$\int_\R |x u_0(x)|dx<\infty\qquad\text{or}\qquad
          \int_\R|xu_0'(x)|dx<\infty,
$$
        then there exist $C,\,N>0$, only depending on
        $(b,\sigma,T ,u_0)$, such that
        \[
            \norm{P_tu_0-\bar P_{t/n}^n u_0}_1\leq \frac{Ct}{n}
                \quad \text{for all}\quad n\geq N.
        \]
    \end{prop}

      The first part of Proposition \ref{prop lgv rate} is
        consistent with  Theorem \ref{thm disc PI converges}, while
        the second part shows that a linear rate can be obtained under
        additional assumptions.
The rate is consistent with estimates for SDEs with
bounded coefficients, see e.g. Corollary 2.7 in \cite{BallyTalayII1995} and Theorem 2.3 in \cite{GobetLa2008}.

   We give a direct proof of Proposition \ref{prop lgv rate} in
   Appendix \ref{sec rate}. Even though we do not have a counter example, the very
        explicit proof seems to indicate that convergence rates can
        only be obtained under additional assumptions on the densities.
    \begin{rem}
        The Cauchy distribution $u_0(x)=\frac{1}{\pi(1+x^2)}$ has no moments,
        but still satisfies $\int_\R|xu_0'(x)|dx<\infty$
        \footnote{An absolutely continuous $L^1$-functions not satisfying
        either of the conditions:
        $u_0(x)=\sum_{k\neq 0}\chi_{I_k}(x)(1-k^2|x-k|)$,
        where $I_k:=[k-1/k^2, k+1/k^2]$ and $k\in \mathbb Z\setminus\set{0}$.}.
    \end{rem}
\end{exam}

\section{Proofs of Theorems \ref{thm solve FPK}, \ref{thm FPK PDF},
  and \ref{thm disc PI converges}}
\label{sec pf convg of disc pi}
We first prove Theorems \ref{thm solve FPK} and \ref{thm disc PI converges}
 since the latter is needed in the proof of Theorem \ref{thm FPK PDF}.
We will use semi-group theory and we refer to
\cite{engel2006one} for more information and precise definitions of the
concepts used below. A crucial step in the proof is to obtain the new
generation result, Theorem \ref{thm gen res}. From
the generation result we obtain a solution of the FPK equation
\eqref{eq FPK}, and Theorem \ref{thm solve FPK} follows. Then we prove
that FPK operator $A$ is {\em dissipative}. Armed with generation and
dissipativity, we use the Chernoff product formula
to show that the solution of the discrete path integration method
converges to the solution of the FPK equation \eqref{eq FPK} and
hence Theorem \ref{thm disc PI converges} follows.
Finally we show Theorem \ref{thm FPK PDF}, that the solution
of \eqref{eq FPK} is just the PDF of the solution of \eqref{eq_general_sde}.

In the rest of this section we let $A$ be defined in \eqref{eq def of A}
and take $D(A):=C_c^{\infty}(\Real^d)$. First we state the generation result.
\begin{thm}\label{thm gen res}
    Assume (C1) and (C2). Then $(A, D(A))$ is closable and its closure $(\bar{A}, D(\bar{A}))$
    generates a contraction semigroup $P_t$ on $(L^1(\Real^d), \norm{\cdot}_1)$.
\end{thm}
This result is relatively technical, and the proof is left to Section
\ref{sec gen result}.
\begin{lem}\label{thm:strcont}
    Assume (C1) and (C2). Then the semigroup $P_t$ from Theorem
    \ref{thm gen res} is a strongly continuous semigroup on
    $(L^1(\Real^d),\norm{\cdot}_1)$.
\end{lem}
\begin{proof}
    The result follows by Proposition I.5.3 in \cite{engel2006one}
    with $T_t=P_t$ if we can verify that
    there exist $\delta>0$, $M\geq 1$ and
    a dense subset $D\subset X$ such that
    \begin{enumerate}
        \item $\norm{T_t}\leq M$ for all $t\in[0,\delta]$,
            \medskip
        \item $\displaystyle \lim_{t\to 0+}T_tx=x$ for all $x\in D$.
    \end{enumerate}
    Since $P_t$ is a contraction semigroup by Theorem \ref{thm gen res},
    condition $(1)$ is satisfied with $M=1$. Moreover, assumption $(2)$
    holds on the dense set $D=D(A)$,  since by the definition of the generator $A$:
    \begin{align*}
          &\lim_{t\to 0+}\left( P_t x-x\right)  =  \lim_{t\to 0+}t\left(
            P_t x-x\right)/t = \left(\lim_{t\to 0+}t\right)\lim_{t\to
            0+}\frac{1}{t}\left( P_t x-x\right)= 0\cdot Ax = 0.
    \end{align*}
    Hence $(1)$ and $(2)$ holds, and we conclude that $P_t$ is strongly continuous.
\end{proof}

\begin{proof}[Proof of Theorem \ref{thm solve FPK}]
    In view of Theorem \ref{thm gen res} and Lemma \ref{thm:strcont},
    the closure of $(A, D(A))$ generates a strongly continuous semigroup
    $P_t$ on $(L^1(\Real^d), \norm{\cdot}_1)$. Let $u(t,x):=(P_t
    u_0)(x)$, and note that $u\in
    C([0,T];L^1(\R^d))$ by strong continuity. Then by
    e.g. Proposition II.6.4 in
    \cite{engel2006one}, $u$ is the
    unique mild solution of the FPK equation \eqref{eq FPK}.
\end{proof}

To proceed, we introduce the following definition.
\begin{defn}
    A linear operator $(B,D(B))$ on a Banach space $(\mathbb X, \norm{\cdot})$ is  \emph{dissipative} if
    \[
        \norm{(\lambda-B)u}\geq \lambda\norm{u}
    \]
    for all $\lambda>0$ and all $u\in D(B)$.
\end{defn}
One can refer to Section II.3 of \cite{engel2006one} for more discussion of this concept.

\begin{prop}\label{prop dissipative}
    Assume (C1). Then the operator $(A, D(A))$ is dissipative on $(L^1(\Real^d), \norm{\cdot}_1)$.
\end{prop}

\begin{proof}
    Let $u\in D(A)$ and define
    $E^\pm:=\set{x\in \Real^d: \pm u(x)> 0}$.
    Since $E^\pm\subset\R^d$, $|u|=\pm u$ on $E^\pm$,
    and $|(\lambda-A)u|\geq (\lambda-A)(\pm u)$,
    \begin{align*}
        \norm{(\lambda-A)u}_1
         & \geq \int_{E^+\cup E^-}|(\lambda-A)u| \geq \int_{E^+\cup E^-} (\lambda-A)|u|.
    \end{align*}
    We claim that
    \begin{align}\label{claim}
        \int_{E^\pm} A|u|\leq 0,
    \end{align}
    and hence the proposition follows since
    $$        \norm{(\lambda-A)u}_1\geq \int_{E^+\cup E^-} \lambda |u|=\lambda\|u\|_1.$$

    Now we prove  claim \eqref{claim} for the $E^+$ case.
    If $E^+\neq\emptyset$, we can
    approximate it by sets $E^+_{\eps_n}=\set{x:u(x)>\eps_n}$,
    $0<\eps_n\to0$, with $C^1$ boundaries. This can be done since by
    Sard's theorem and the implicit function theorem, $E^+_{\eps}$ has
    $C^1$ boundary for a.e. $0<\eps<\max u^+$. Note that
    $E^\pm_\eps\subset E^\pm\subset \text{supp}\, u$ which is compact.

    Then we write $A$ in divergence form, cf. \eqref{eq def of A},
    \[
         A\phi
         = \frac{1}{2}\,\text{div}\left(a\nabla \phi
            +(\text{div }a_1,\cdots,\text{div }a_d)^T\phi-2b\phi \right),
    \]
    and use the divergence theorem. We consider first the diffusion term,
    \[
       \int_{E^+_\eps}\text{div }(a\nabla u) = \int_{\partial
         {E^+_\eps}} a\nabla u \cdot\mathbf{n},
    \]
    where $\mathbf{n}$ the exterior unit normal vector of $\partial E^+_\eps$.
    Since $\partial {E^+_\eps}$ is an $\eps$-level set of $u$ and $u$
    is decreasing in the outward direction at $\partial {E^+_\eps}$,
    $\nabla u=-\beta \mathbf{n}$ for some $\beta \geq 0$.
    Then since $a=\sigma\sigma^T$ is positive semi-definite
    \[
         a\nabla u\cdot \mathbf{n}=-\beta \mathbf{n}^Ta\mathbf{n}\leq 0,
    \]
    and hence
    \[
        \int_{E^+_\eps}\text{div }(a\nabla u)\leq 0.
    \]
    Next we  estimate the convection part. Since $u=\eps$ on $\partial
    E^+_\eps$,
    \begin{align*}
       &\int_{E^+_\eps}\text{div}\Big((\text{div\,}a_1,\cdots,\text{div\,}a_d)^Tu-2bu\Big) =  \int_{\partial
            E^+_\eps}((\text{div}\,a_1,\cdots,\text{div\,}a_d)^T-2b)
          u\cdot \mathbf{n} \\
        & = \eps \int_{\partial E^+_\eps}((\text{div
        }a_1,\cdots,\text{div }a_d)^T-2b) \cdot\mathbf{n} = \eps
        \int_{E^+_\eps}\text{div}\Big((\text{div }a_1,\cdots,\text{div
        }a_d)^T-2b\Big),
    \end{align*}
    which is bounded by $C \eps$ for some $C=C(a,b,\text{supp} u)$ since $a$ and $b$
    are smooth and $u$ has compact support. Hence
    \begin{equation}\label{eq div est}
      \int_{E^+_\eps}A u \leq C\eps.
    \end{equation}

    Note that
    $\chi_{E^+_\eps}\to \chi_{E^+}$ a.e. Since  $Au$ belongs to
    $L^1(\Real^d)$, we may then use the dominated convergence theorem
    to pass to the limit in \eqref{eq div est} as $\eps\to 0^+$  and obtain
    \eqref{claim}. The $E^-$
    case is similar and will be omitted.  The proof is complete.
\end{proof}

Now we proceed to prove Theorem \ref{thm disc PI converges}. To do
that we need the \emph{Chernoff Product Formula}
(e.g. Theorem III.5.2 in \cite{engel2006one}):
\begin{thm}\label{thm chernoff prod formula}
    Let $\mathbb X$ be a Banach space, $V: \Real_+ \to \mathcal{L}(\mathbb X)$ a
    function such that $V(0)=id$ and $\|V^n(t)\|\leq M$ for all $t\geq 0$,
    all $n \in\mathbb{N}$, and some $M\geq 1$.
    Assume
    \begin{equation}\label{eq lim of V(t)}
      Bu:=\lim_{t\to 0+}\frac{V(t)u-u}{t}
    \end{equation}
    exists for all $u\in D\subset \mathbb X$, where $D$ and $(\lambda-B)D$
    are dense subspaces in $\mathbb X$ for some $\lambda>0$.
    Then the closure $\bar{B}$ of $B$ generates a bounded strongly continuous
    semigroup $\set{T(t): t\geq 0}$, given by
    \begin{equation}\label{eq group from V(t)}
      T(t)u=\lim_{n\to +\infty}V^n\left(t/n\right) u,
    \end{equation}
    for all $u\in \mathbb X$ uniformly for $t\in[0,T]$.
\end{thm}

This result does not require the approximation operators $V(t)$ to be
a semigroup. This is important since the discrete path integration operators
$\set{\bar{P}_{\tau}}$ are bounded on the Banach space $L^1(\Real^d)$
by Proposition \ref{prop disc PI contract}.
But in contrast to their continuous counterpart, they no longer
form a semigroup of operators\footnote{E.g. in Example \ref{exam_lgv}, $\bar{P}^2_{\tau}\neq
\bar{P}_{2\tau}$ since
\begin{align*}
  (\bar{P}^2_{\tau}f)(y) & = \int_{\Real}\frac{f(x)}{\sqrt{4\pi\tau \sigma^2(1+b\tau+b^2\tau^2/2)}}
        \exp\set{-\frac{(y-(1+2b\tau+b^2\tau^2)x)^2}{4\tau \sigma^2(1+b\tau+b^2\tau^2/2)}}dx, \\
  (\bar{P}_{2\tau}f)(y) & =\int_{\Real}\frac{f(x)}{\sqrt{4\pi\tau \sigma^2}}
        \exp\set{-\frac{(y-(1+2b\tau)x)^2}{4\tau \sigma^2}}dx.
\end{align*}}.

\begin{proof}[Proof of Theorem \ref{thm disc PI converges}]
    We use Theorem \ref{thm chernoff prod formula}
    with $\mathbb X=L^1(\Real^d)$, $V(h)=\bar P_h$, $B=A$, and $D=D(A)$. Let us verify
    the conditions. By Proposition \ref{prop lim of disc PI},
    $$A u = \lim_{h\to 0} \frac{\bar P_h u-u}h$$
    for all $u\in D(A)$. By the definition of $D(A)$ and Theorems
    \ref{prop dissipative} and \ref{thm gen res}, $A$ is densely defined, dissipative,
    and the closure $\bar A$ generates a contraction
    semigroup $P_t$ on $L^1(\Real^d)$. Then by the Lumer-Phillips Theorem
    (e.g. Theorem II.3.15 in \cite{engel2006one}), this implies that
    $(\lambda-A)D(A)$ is dense in $L^1(\R^d)$ for some
    $\lambda>0$.

    Since the semigroup is strongly continuous by Lemma
    \ref{thm:strcont}, we have verified all the conditions for Theorem
    \ref{thm chernoff prod formula}. Hence hence can conclude that
    $$\norm{u(t,\cdot)-\bar{u}_n}_1=\norm{P_t u_0-\bar{P}_{t/n}^n
      u_0}_1\to 0,$$
    as $n\to \infty$ uniformly in $[0,T]$. The proof is complete.
\end{proof}

Now we prove Theorem \ref{thm FPK PDF}, using among other things
Theorem \ref{thm disc PI converges}.
\begin{proof}[Proof of Theorem \ref{thm FPK PDF}]
    We will show that \eqref{eq def of PDF} holds.

    Assume first that $\phi\in C_b^3(\Real^d)$.
    In view of the Euler-Maruyama scheme \eqref{eq EM disc}
    corresponding to $Y_t$ with
    $\tau=\frac t n$, the PDF of its solution $X_{n}$ is
    $\bar{u}_n:=\bar{P}^n_{t/n}u_0$, and hence
    \begin{equation}\label{eq disc pdf}
        \Epc\phi(X_n)=\int_{\Real^d}\phi(x)
        \bar{u}_n(x)dx\qquad\text{for any}\qquad\phi\in C_b^0(\R^d).
    \end{equation}
    The right hand side is close to what we want, since by Theorem
    \ref{thm disc PI converges},
    \[
        \abs{\int \phi \,\bar{u}_n-\int \phi \, u}
            \leq \norm{\phi}_{\infty}\sup_{t\in
              [0,T]}\norm{u-\bar{u}_n}_1\to 0\qquad
            \text{as}\qquad n\to \infty.
    \]

    To see that the left hand side of \eqref{eq disc pdf} is also close to
    the sought after expression, we will use the weak convergence of the
    Euler-Maruyama method. Since our initial distribution $u_0$ is not
    assumed to have second moments, we need to introduce the following
    auxiliary SDE and Euler-Maruyama scheme for a fixed $x\in \Real$:
    \[
        \begin{cases}
            d \tilde Y_t = b(\tilde Y_t)\, dt + \sigma(\tilde Y_t) \, dB_t,\\
            \tilde Y_0 = x,
        \end{cases}\quad\text{and}\qquad
        \begin{cases}
            \tilde X_{n+1} = \tilde X_n + b(\tilde X_n) \tau +
                \sigma(\tilde X_n)\Delta B_{n},\\
            \tilde X_0=x.
        \end{cases}\!\!\!\!\!\!\!\!\!\!
    \]
    Both problems have unique (strong) solutions (Proposition
    \ref{prop SDE exist and uniq}), so we may define
    \[
        v(t,x):=\Epc\phi(\tilde{Y}_t)\qquad\text{and}\qquad v_n(x):=\Epc\phi(\tilde{X}_n).
    \]
    By the strong Markov property of solutions of SDEs, cf.
    Theorem V.32 in \cite{protter2005stochastic},
    \[
        \Epc\phi(Y_t)=\Epc[\Epc(\phi(Y_t)|Y_0=Z)]=\Epc[v(t,Z)].
    \]
    Similarly, since $\set{X_n}$ is a Markov chain,
    \[
        \Epc\phi(X_n)=\Epc[\Epc(\phi(X_n)|X_0=Z)]=\Epc[v_n(Z)].
    \]
    Therefore
    \[
        \abs{\Epc\phi(Y_t)- \Epc\phi(X_n)}
            \leq \Epc\abs{v(t,Z)-v_n(Z)}=\int_{\Real^d}\abs{v(t,x)-v_n(x)}u_0(x)dx.
    \]
    The latter integrand is dominated by the integrable function
    $2\norm{\phi}_{\infty}$, and by the weak convergence of the Euler-Maruyama scheme,
    cf. Theorem 14.1.5 in \cite{kloeden2013numerical},
    \[
        \abs{v(t,x)-v_n(x)}\to 0\qquad \text{as}\qquad n\to \infty\qquad\text{for every $x\in\R^d$}.
    \]
    Hence by (C3) and the dominated convergence theorem, it follows that
    \[
           \int_{\Real^d}\abs{v(t,x)-v_n(x)}u_0(x)dx
            \to 0\qquad \text{as}\qquad n\to \infty.
    \]
    Note that the Theorem of \cite{kloeden2013numerical} requires the initial
    probability measure to have finite second moments which is trivially met for the
    auxiliary problems above.

    From all of the estimates above, we can the conclude that
    \begin{align*}
        \abs{\Epc\phi(Y_t)- \int \phi \, u}
         &\leq \abs{\Epc\phi(Y_t)- \Epc\phi(X_n)}
            + \abs{\int \phi \, \bar{u}_n-\int \phi  \,u} \to 0\quad
            \text{as}\quad n\to \infty,
    \end{align*}
    and hence that
    \[
        \Epc\phi(Y_t)=\int_{\Real^d} \phi(x) u(t,x)dx
            \qquad \text{for all}\qquad \phi\in C_b^3(\Real^d).
    \]

    Now we assume $\phi\in C_b^0(\Real^d)$ and conclude the proof by a
    regularization argument.
    Let $\rho_{\eps} \in C_c^{\infty}$, $\eps>0$, be the standard mollifier and
    $\phi_{\eps}:=\rho_{\eps}*\phi$. Then by e.g. Section 7.2 in
    \cite{gilbarg2015elliptic}, $\phi_{\eps}\in C_b^\infty(\Real^d)$,
    $\norm{\phi_{\eps}}_{\infty}\leq\norm{\phi}_{\infty}$, and
    $\phi_{\eps}\to \phi$ a.e. Observe that
    \begin{align*}
       & \abs{\int_{\Real^d}\phi(x)u(t,x)dx-\Epc\phi(Y_t)} \\
       & \leq \abs{\int \phi\, u- \int \phi_{\eps}u}+\abs{\int \phi_{\eps}u-\Epc\phi_{\eps}(Y_t)}
            +\abs{\Epc\phi_{\eps}(Y_t)-\Epc\phi(Y_t)}.
    \end{align*}
    The second term is zero by our result for $\phi \in C_b^3(\Real^d)$,
    and the two remaining terms tend to zero as $\eps\to 0+$ by the
    dominated convergence theorem. The proof is complete.
\end{proof}

\section{Proof of Proposition \ref{prop lim of disc PI}}\label{Pf2}
This section is devoted to the proof of Proposition \ref{prop lim of disc PI}.
    The key idea comes from the
    proof of Lemma 16 in \cite{Butko2016}, and we extend it to the case
    with unbounded coefficients.

\begin{proof}
    Since $u\in C_c^\infty(\Rd)$,
    there is $R>0$ such that $\supp u \subset B_R(0)$.
    The PDF after one-step Euler-Maruyama scheme is
    $\bar P_\tau u$. By Lemma 2.7,
    \begin{align*}
       & \mathcal F \left(\bar P_\tau u - u - \tau Au \right) \\
       & = \int_{\Rd}u(x)e^{i\xi x}
            \left(\hat{\bar k}(\xi, x,\tau) -1 - \tau(ib(x)\xi-\frac 12
                \xi^T a(x)\xi)\right) dx\\
       & = \tau^2 \iint_0^1 u(x) e^{i\xi x} (1-\theta)\left(ib(x)\xi -\frac 12
            \xi^T a(x)\xi\right)^2 e^{\theta\tau\left(ib(x)\xi -\frac 12
            \xi^T a(x)\xi\right) } d\theta dx.
    \end{align*}
    Without loss of generality, we assume $d=1$.
    Note that $D_x^k e^{i\xi x}= i^k \xi^k e^{i\xi x}$, we have
    \begin{align*}
       & \mathcal F \left(\bar P_\tau u - u - \tau Au \right) \\
       & = \tau^2\int_0^1 (1-\theta)\int_\R e^{i\xi x}
            \left[ D_x^2\left( b(x)u(x)e^{\theta\tau\left(ib(x)\xi -\frac 12
            a(x)\xi^2\right)}\right)\right.\\
         &   \qquad + D_x^3\left( a(x)b(x)u(x)
                e^{\theta\tau\left(ib(x)\xi -\frac 12
            a(x)\xi^2\right)} \right) \\
          &  \qquad \left. + \frac 14 D_x^4\left( a^2(x)u(x)
                e^{\theta\tau\left(ib(x)\xi -\frac 12
            a(x)\xi^2\right)} \right)  \right] dx d\theta\\
         & = \tau^2\int_0^1 (1-\theta)\int_\R e^{i\xi x}
            e^{\theta\tau\left(ib(x)\xi -\frac 12 a(x)\xi^2\right)}
            \sum_{m,n=0}^4 \psi_{m,n}(x) (i\theta t \xi)^m (-\theta t \xi^2)^n
                dx d\theta,
    \end{align*}
    where $\psi_{m,n}\in C_c(\R)$ and for all $x\in \R$,
    \[
        \abs{\psi_{m,n}(x)}\leq C_{m,n} \sum_{k=0}^4 \abs{D^k u(x)},
    \]
    and the positive constant $C_{m,n}$ depends on $b$, $\sigma$,
     their derivatives up to fourth order, and $R$.

    Now we perform inverse Fourier transform and estimate it in $L^1$-norm
    \begin{align*}
       & \norm{\bar P_\tau u - u - \tau Au}_1 \\
       & = \tau^2 \int_\R \Big | \int_\R \int_0^1 (1-\theta)\int_\R e^{i\xi x}
            e^{\theta\tau\left(ib(x)\xi -\frac 12 a(x)\xi^2\right)}\\
       &     \qquad \cdot \sum_{m,n=0}^4 \psi_{m,n}(x)
            (i\theta t \xi)^m (-\theta t \xi^2)^n dx d\theta d\xi \Big | dy\\
       &  \leq \tau^2 \int_0^1 \int_\R\Big | \int_\R \int_\R e^{i\xi (x-y)}
            e^{\theta\tau\left(ib(x)\xi -\frac 12 a(x)\xi^2\right)}\\
       &     \qquad \cdot \sum_{m,n=0}^4 \psi_{m,n}(x)
            (i\theta t \xi)^m (-\theta t \xi^2)^n dx d\xi \Big | dy d\theta \\
       &   = \tau^2 \int_0^1 \int_\R\Big | \int_\R \int_\R
                e^{i\eta \frac{x-y}{\sqrt{\theta \tau}}}
                e^{ib(x)\sqrt{\theta \tau}\eta -\frac 12 a(x)\eta^2}\\
       &     \qquad \cdot \sum_{m,n=0}^4 \psi_{m,n}(x)
            (i\sqrt{\theta \tau}\eta)^m (-\eta^2)^n
                \frac{dx}{\sqrt{\theta \tau}}d\eta \Big | dy d\theta \\
       &   = \tau^2 \int_0^1 \int_\R\Big | \int_\R \int_\R
                e^{i\eta z}
                e^{ib(y+ \sqrt{\theta \tau}z)\sqrt{\theta \tau}\eta -\frac 12 a(y+ \sqrt{\theta \tau}z)\eta^2}\\
       &     \qquad \cdot \sum_{m,n=0}^4 \psi_{m,n}(y+ \sqrt{\theta \tau}z)
            (i\sqrt{\theta \tau}\eta)^m (-\eta^2)^n
                 dz d\eta \Big | dy d\theta.
    \end{align*}
    After the above changes of variables, we continue the estimate
    \begin{align*}
       & \norm{\bar P_\tau u - u - \tau Au}_1 \\
       &   \leq \tau^2 \int_0^1 \int_\R\Big | \int_\R
                \sum_{m,n=0}^4 \psi_{m,n}(y+ \sqrt{\theta \tau}z)
                (i\sqrt{\theta \tau})^m (-1)^n\\
       &     \qquad \cdot \left(\int_\R e^{i\eta z} \eta^{m+2n}
            e^{ib(y+ \sqrt{\theta \tau}z)\sqrt{\theta \tau}\eta
                    -\frac 12 a(y+ \sqrt{\theta \tau}z)\eta^2}
            d\eta \right) dz\Big | dy d\theta\\
       &   = \tau^2 \int_0^1 \int_\R\Big | \int_\R
                \sum_{m,n=0}^4 \psi_{m,n}(y+ \sqrt{\theta \tau}z)
                (\sqrt{\theta \tau})^m\\
       &     \qquad \cdot \left(D_v^{m+2n} \int_\R
            e^{i(v+ b(y+ \sqrt{\theta \tau}z)\sqrt{\theta \tau})\eta}
            e^{-\frac 12 a(y+ \sqrt{\theta \tau}z)\eta^2}
            d\eta \right)_{v=z} dz\Big | dy d\theta\\
       &   = \tau^2 \int_0^1 \int_\R\Big | \int_\R
                \sum_{m,n=0}^4 \psi_{m,n}(y+ \sqrt{\theta \tau}z)
                (\sqrt{\theta \tau})^m\\
       &     \qquad \cdot \left(D_v^{m+2n}
            \frac{1}{\sqrt{2\pi a(y+ \sqrt{\theta \tau}z)}}
            e^{-\frac{(v+b(y+ \sqrt{\theta \tau}z)\sqrt{\theta \tau})^2}
                {2a(y+ \sqrt{\theta \tau}z)}}
            \right)_{v=z} dz\Big | dy d\theta.
    \end{align*}

    Now there are no more oscillatory integrals, and therefore
    we can continue the estimate as follows
    \begin{align*}
       & \norm{\bar P_\tau u - u - \tau Au}_1 \\
       & \leq \tau^2 \sum_{m,n=0}^4 \int_0^1 \int_\R \int_\R
        \abs{\psi_{m,n}(y+ \sqrt{\theta \tau}z)}\\
        &     \qquad \cdot \abs{ \left(D_v^{m+2n}
            \frac{1}{\sqrt{2\pi a(y+ \sqrt{\theta \tau}z)}}
            e^{-\frac{(v+b(y+ \sqrt{\theta \tau}z)\sqrt{\theta \tau})^2}
                {2a(y+ \sqrt{\theta \tau}z)}}
            \right)_{v=z}} dy dz d\theta\\
       & = \tau^2 \sum_{m,n=0}^4 \int_0^1 \int_\R \int_\R
        \abs{\psi_{m,n}(x)  \left(D_v^{m+2n}
            \frac{1}{\sqrt{2\pi a(x)}}
            e^{-\frac{(v+b(x)\sqrt{\theta \tau})^2}
                {2a(x)}}
            \right)_{v=z}} dx dz d\theta\\
       & \leq \tau^2\sum_{m,n=0}^4  \int_0^1  \int_\R \int_\R
        \abs{\psi_{m,n}(x)} \\
        & \qquad \cdot \frac {C_{m,n}' }{\sqrt{2\pi a(x)}}   \left( 1+
            \frac{(z+\sqrt{\theta \tau}b(x))^{2(m+2n)}}{a^{m+2n}(x)}
            \right) e^{-\frac{(z+\sqrt{\theta \tau}b(x))^2}{2a(x)}}
                dx dz d\theta\\
        & \leq C_1 \tau^2\sum_{m,n=0}^4  \int_0^1
         \int_\R  \abs{\psi_{m,n}(x)}\\
        & \qquad \cdot \left( \int_\R \frac {1}{\sqrt{2\pi a(x)}}
            \left( 1+  \left(\frac{z+\sqrt{\theta \tau}b(x)}{\sqrt{a(x)}}\right)^{2(m+2n)}\right)
            e^{-\frac{(z+\sqrt{\theta \tau}b(x))^2}{2a(x)}}dz\right)
             dx d\theta.
    \end{align*}
    Finally we make a further change of variable
    $\tilde z= \frac{z+\sqrt{\theta \tau}b(x)}{\sqrt{a(x)}}$ and obtain
    \begin{align*}
       & \norm{\bar P_\tau u - u - \tau Au}_1 \\
        & \leq C_1 \tau^2\sum_{m,n=0}^4  \int_0^1
         \int_\R  \abs{\psi_{m,n}(x)} \left( \int_\R \frac {1}{\sqrt{2\pi }}
            \left( 1+  \tilde z^{2(m+2n)}\right)
            e^{-\frac{\tilde z^2}{2}}d\tilde z\right)
             dx d\theta\\
        & \leq C_2 \tau^2 \sum_{m,n=0}^4 \norm{\psi_{m,n}}_1\\
        & \leq C \tau^2 \sum_{k=0}^4 \norm{D^k u}_1.
    \end{align*}
    That is, $\norm{\frac{\bar{P}_{\tau}u-u}{\tau}-Au}_1
            \leq C(b,\sigma,u) \tau$.
    The proof is complete.
\end{proof}
\begin{rem}
    From the above proof we see that the consistency
    estimate is uniform in the elliptic number $\alpha>0$ in condition (C2).
\end{rem}

\section{Proof of Theorem \ref{thm gen res}}\label{sec gen result}
To prove Theorem \ref{thm gen res} we need first to show that operator
$A$ and its adjoint satisfy certain properties, including being dissipative.

\begin{prop}\label{prop dissipative A}
    Assume (C1). Then the operators $A$ on $D(A):=C_c^{\infty}(\Real^d)$
    is densely defined, dissipative, and closable (possessing a closed extension) in $L^1(\Real^d)$.
\end{prop}
\begin{rem}
    In view of Theorem \ref{thm chernoff prod formula} and Proposition
    \ref{prop dissipative} we only need $A$ to be closable.
    Alternatively we could have shown directly that $A$ is closed with
    respect to the operator norm $\norm{u}_A:=\norm{u}_1+\norm{Au}_1$ as in
    Section 2 of \cite{Fornaro2007747}.
\end{rem}
\begin{proof}[Proof of Proposition \ref{prop dissipative A}]
    The operator $A$ is densely defined since $D(A)$ is dense in
    $L^1(\Real^d)$,  and it is dissipative by Proposition \ref{prop dissipative}.
    Hence it is closable by Proposition II.3.14 in \cite{engel2006one}.
\end{proof}

On the dual space of $(L^1(\Real^d), \norm{\cdot}_1)$,
$(L^{\infty}(\Real^d), \norm{\cdot}_{\infty})$, the adjoint operator
of $A$ with domain $D(A)=C_c^{\infty}(\Real^d)$ is given by the
definition below.
\begin{defn}\label{def A^*}
    The adjoint of $A$ is an operator $A^*: D(A^*)\to
    L^{\infty}(\Real^d)$ defined as
    \begin{equation}\label{eq int def of A^*}
        \int_{\Real^d} u A^*f=\int_{\Real^d} fAu \text{\qquad for all\qquad }f\in\ D(A^*) \text{ and }u\in D(A),
    \end{equation}
    where $D(A^*)=\set{f\in L^{\infty}: \exists \,g \in L^{\infty}\text{ such that }
        \int fAu=\int gu \text{ for all }u\in D(A)}$.
\end{defn}

\begin{rem}\label{rem A^*}
    By integration by parts in \eqref{eq def of A},
    $A^* f$ is given by \eqref{eq adj op on Cb} for any $f\in C_b^2(\Real^d)$.
    Moreover $f\in D(A^*)$ if and only if there exists a $g\in
    L^{\infty}(\Real^d)$ such that
    \begin{equation}\label{eq A^*f=g}
        A^* f=g \qquad\text{in}\quad\mathcal D' \text{ (in distributions)}.
    \end{equation}
\end{rem}

\begin{prop}\label{prop dissipative A adjoint}
    Assume (C1) and (C2). Then the adjoint operator $(A^*, D(A^*))$ is
    dissipative on $(L^{\infty}(\Real^d),\norm{\cdot}_{\infty})$.
\end{prop}

The proof is quite technical, so we postpone it to the end of this section.
Now we can employ the following result to prove that $A$ generates a contraction semigroup.

\begin{proof}[Proof of Theorem \ref{thm gen res}]
    Note that $A$ with $D(A)=C_c^\infty(\R^d)$ is a bounded densely defined
    operator on $L^1(\R^d)$, and that both $A$ and its adjoint $A^*$ are
    dissipative by Propositions \ref{prop dissipative A} and
    \ref{prop dissipative A adjoint}. Hence the closure $\bar A$ generates
    a contraction semigroup on
    $(L^1(\Real^d),\norm{\cdot}_1)$ by Corollary II.3.17 in
    \cite{engel2006one}. The proof is complete.
\end{proof}

In the rest of this section we prove Proposition \ref{prop dissipative A adjoint}.

\begin{lem}\label{lem ell reg}
    Assume (C1) and (C2), then
    \[
        D(A^*)\subset \set{f\in C^1(\Real^d)\cap L^{\infty}(\Real^d) :
            \nabla^2 f\in\bigcap_{1<p<\infty}L^p_{loc}(\Real^d)}.
    \]
\end{lem}
\begin{rem}
    Under different assumptions, similar characterisations of the
    maximal domain can be found in Chapter 3 of \cite{lunardi2012analytic}
    and Section 3 of \cite{Fornaro2007747}.
\end{rem}
\begin{proof}
    For any $f\in D(A^*)$, there exists $g\in L^{\infty}(\Real^d)$
    such that \eqref{eq A^*f=g} holds by Remark \ref{rem A^*}.
    By the regularity of the coefficients (see (C1)), we rewrite \eqref{eq int def of A^*} as
    \begin{align*}
      \notag & \int f \sum_{i,j=1}^d \partial_i \left(\frac 12
        a_{ij} \partial_j u\right)  \\
      & =\int \left(g+f \sum_{i=1}^d \partial_i b_i -\frac 12 f \sum_{i,j=1}^d \partial_i\partial_ja_{ij}\right)u
            +  f\sum_{i=1}^d \left( b_i-\frac 12  \sum_{j=1}^d \partial_j a_{ij}\right) \partial_i u,
    \end{align*}
    On the right hand side, the coefficients of $u$ and $\partial_i u$ belong to
    $L^p_{loc}(\Real^d)$ for all $1<p<\infty$.
    We can therefore apply Theorem 1.5 in \cite{ZhangBao2013} to show
    that $f\in W^{1,p}_{loc}(\Real^d)$ for all $1<p<\infty$.
    With this extra regularity and integration by parts,
    \[
        \int f  \sum_{i,j=1}^d \partial_i \left(\frac 12 a_{ij}\partial_j u\right)
            =\int  \left(g+\frac 12\sum_{i,j=1}^d \partial_i f\partial_ja_{ij}-\sum_{i=1}^d b_i \partial_i f \right)u.
    \]
On the right hand side, the coefficients of $u$ belong to
    $L^p_{loc}(\Real^d)$ for all $1<p<\infty$ while the coefficient of
    $\partial_i u$ is zero. By Proposition 1.1 in \cite{ZhangBao2013}
    we then have that $f\in W^{2,p}_{loc}(\Real^d)$ for all $1<p<\infty$.

    Finally, by taking $p>d$ in the Sobolev embedding theorem, see
    e.g. Theorem 7.10 in
    \cite{gilbarg2015elliptic}, we find that $\nabla f\in
    C(\Real^d)$. The proof is complete.
\end{proof}

The next lemma is crucial.
\begin{lem}\label{lem critical seq of bounded function}
    Assume $f\in D(A^*)$ and  $m:=\underset{\Real^d}{\sup}\, f<\infty$.
    Then there exists a sequence $\set{x_n}\subset \Real^d$ such that
    \[
      \lim_n f(x_n)= m,\quad\lim_n (1+\abs{x_n})\abs{\nabla
        f(x_n)}=0,\quad\text{and}\quad \lim_n \frac{|x_n|^2}{n} =0.
    \]
    Moreover, $\nabla^2 f(x_n)$ exists and
    \[
      \nabla^2 f(x_n) \leq \frac{3}{n}I_d,
    \]
    where $I_d$ is the identity matrix on $\Real^d$.
\end{lem}

Since the second derivatives are only defined almost everywhere,
we need the so-called Bony maximum principle
\cite{Bony1967, Lions1983} to prove the above
lemma. We quote the version given in Proposition 1.2.12
in \cite{garroni2002second}.
\begin{lem}\label{prop max principle}
    Suppose $p\geq d$ and a function $w\in W^{2,p}_{loc}(\Real^d)$
    achieves a local maximum at point $x_0$,
    then for all $\xi\in \Real^d$,
    \[
        \lim_{r\to 0} \left( \underset{B_r(x_0)}{\essinf}
            \sum_{i,j=1}^d \xi_i (\partial_i\partial_j w) \xi_j \right)\leq 0,
    \]
where $B_r(x_0)$ denotes the ball of radius $r$ centered at $x_0$.
\end{lem}

The key idea of the following proof now comes from Lemma 2.3 in
    \cite{Jakobsen2002}.

\begin{proof}[Proof of Lemma \ref{lem critical seq of bounded function}]
    For each $n\in\mathbb{N}$, we define
    \[
        g_n(x):= f(x)-\frac{|x|^2}{n}.
    \]
    Hence $g_n$ is bounded from above and tends to $-\infty$
    when $|x|\to \infty$. Therefore there exists $x_n'$ for each $g_n$
    such that
    \[
        g_n(x_n')=\max_{\Real^d}g_n=:m_n.
    \]
    Now for arbitrary $\epsilon >0$, we can find first $x_{\epsilon}\in \Real^d$
    and  then $N\in \N$ such that
    \[
       f(x_{\epsilon})>m-\epsilon\qquad\text{and}\qquad
        \frac{\abs{x_{\epsilon}}^2}{n}<\epsilon\;\;\text{ for all }n>N.
    \]
    Hence for all  $n>N$,
    \begin{align*}
      m & \geq m_n=f(x_n')-\frac{|x_n'|^2}{n}=\max_{\Real^d}g_n \\
       & \geq g_n(x_{\epsilon})=f(x_{\epsilon})-\frac{\abs{x_{\epsilon}}^2}{n}\\
       & \geq m-2\epsilon.
    \end{align*}
    Therefore $m_n\to m$ as $n\to \infty$, and hence $f(x_n')\to m$
    and $\frac{|x_n'|^2}{n}\to 0$ as $n\to \infty$.

    Next we analyse the derivatives. By definitions of $x_n'$ and $g_n$,
    \[
        0=\nabla g_n(x_n')=\nabla f(x_n')-\frac{2}{n}x_n'.
    \]
    Hence, as $n\to \infty$,
    \[
        (1+\abs{x_n'})\abs{\nabla f(x_n')}
            =(1+\abs{x_n'})\frac{2}{n}|x_n'|\leq \frac{1}{n}\left( 1+3|x_n'|^2 \right)\to 0.
    \]

    Finally for the second derivatives,
    \[
        \nabla^2 g_n(x)=\nabla^2 f(x)-\frac{2}{n}I_d,\text{ a.e.}
    \]
    Since $x_n'$ is maximum point for $g_n$, by Proposition
    \ref{prop max principle},
    \[
        \lim_{r\to 0} \left(\underset{B_r(x_n')}{\essinf}
         \, \nabla^2 g_n \right)\leq0,  \qquad\text{and hence}\qquad
       \lim_{r\to 0} \left(\underset{B_r(x_n')}{\essinf}
         \,\nabla^2 f\right) \leq \frac{2}{n}I_d.
    \]

    By the definition of the infimum and regularity of $f$, we
    can find another sequence $\set{x_n}$ such that $\nabla^2 f(x_n)$ is defined,
    $|x_n-x_n'|<1/n$, and
    \[
        \nabla^2 f(x_n) \leq \frac{3}{n}I_d.
    \]
    Using that $f\in C^1$, we can take a further subsequence such that also
    \[
        |f(x_n)-f(x_n')|+|\nabla f(x_n)-\nabla f(x_n')|< \frac 1n.
    \]

    Combining all the above estimates, we can conclude that
    \begin{align*}
       \lim_n f(x_n) = m, \quad  \lim_n (1+\abs{x_n})\abs{\nabla
         f(x_n)}=0,\quad  \lim_n \frac{|x_n|^2}{n}  =0,
    \end{align*}
    and
    $ \nabla^2 f(x_n)\leq \frac{3}{n}I_d$.
    The proof is complete.
\end{proof}

\begin{proof}[Proof of Proposition \ref{prop dissipative A adjoint}]
    Let $f\in D(A^*)$. We may assume that $m:=\norm{f}_{\infty}=\sup_{\Real^d}f$.
    The case $m=-\inf_{\Real^d}f$ follows in a similar way and is omitted.

    Note that by \eqref{eq adj op on Cb},
    $a_{ij}(x)=\sum^n_{k=1}\sigma_{ik}(x)\sigma_{jk}$
    and $\sigma=(\sigma_1,\cdots,\sigma_n)$, and then
    \[
        A^*f(x)=\sum_{i=1}^d b_i(x)\partial_if(x)+\frac{1}{2}\sum_{k=1}^n \sigma_k^T(x) \nabla^2f(x)\sigma_k(x).
    \]
    Let $\set{x_n}$ be the sequence corresponding to $f$ given by
    Lemma \ref{lem critical seq of bounded function}.
    From the proof it follows that for any fixed $\lambda>0$, we can always take
    \[
        \set{x_n}\subset \set{x\in\Real^d: |\lambda f(x)- A^*f(x)|\leq \norm{\lambda f-A^*f}_{\infty}},
    \]
    because the complement of the latter set has zero Lebesgue measure in $\Real^d$.
    By (C1) and Lemma \ref{lem critical seq of bounded function},
    \[
      \lim_n f(x_n)=m,\quad\lim_n \abs{b_i(x_n)\partial_if(x_n)}
        \leq \lim_n K(1+\abs{x_n})\abs{\nabla f(x_n)}=0, \\
    \]
    and
    \[
      \limsup_n  \sigma_k^T(x_n)  \nabla^2f(x_n)\sigma_k(x_n)
            \leq \lim_n \frac{3}{n}\abs{\sigma_k(x_n)}^2  \leq \lim_n \frac{3}{n}K(1+\abs{x_n}^2)=0.
    \]
    Hence,
    \[
        \limsup_n (-A^*f(x_n))  \geq 0,
    \]
    and then
    \[
      \lambda \norm{f}_{\infty}=\lambda\lim_n f(x_n)
      \leq \limsup_n (\lambda f(x_n)-A^*f(x_n)) \leq \norm{(\lambda-A^*)f}_{\infty}.
    \]
    Since $f\in D(A^*)$ and $\lambda>0$ are both arbitrary, $A^*$ is a dissipative operator,
    and the proof is complete.
\end{proof}

\appendix

\section{Weak Convergence and Convergence of PDFs}
\label{app:weak}
Now we demonstrate the relationship between
 $L^1$-convergence of PDFs and  weak convergence
 of processes.
\begin{prop}\label{L1weak}
    Assume $X_n$ and $X$ are random variables with
    PDFs $u_n$ and $u$ satisfying
    $$\lim_{n\to\infty}\norm{u_n-u}_1=0.$$
    Then $X_n$ converges weakly to $X$,
    \[
       \lim_{n\to\infty} \abs{\Epc \phi(X_n)-\Epc \phi(X)}= 0
        \quad \text{ for all}\quad \phi\in C_b^0(\Real^d).
    \]
\end{prop}
\begin{proof}
    $\abs{\Epc \phi(X_{n})-\Epc \phi(X)}=\abs{\int_{\R^d}\phi(x)(u_n-u)(x)\,dx}
        \leq \norm{\phi}_{\infty}\norm{\bar{u}_n -u}_1$.
\end{proof}

The opposite is not true in general as the following example shows.
\begin{exam}\label{exweak}
    For any $n\in\N$, we define
    \[
        v_n(x)=\big(1+\sin(2n\pi x)\big)\chi_{[0,1]}(x).
    \]
    Evidently, $v_n\geq 0$ and $\int_{\Real}v_n=1$, and
    hence $v_n$ is a PDF.

    For any $n\in\N$, we take $X_n$ to be a random variable with PDF
    $v_n$.
    By the Riemann-Lebesgue Lemma, $X_n$ converges weakly to a random
    variable $X$ with PDF
    $$v(x)\equiv \chi_{[0,1]}(x).$$
    However, $v_n$ does not converge to $v$ in $L^1$ since
    \[
        \norm{v_n-v}_1:=\int_{\Real}\abs{v_n-v}
            =\int_0^1\abs{\sin(2n\pi x)}dx\equiv \frac{2}{\pi}> 0.
    \]
\end{exam}

\section{Proof of Proposition \ref{prop lgv rate}}\label{sec rate}
    Fix $t\in [0, T]$ and define $\tau=t/n$ for $n>|bt|$.
    Without loss of generality,
    we assume that $b\neq 0$.
    We start with deriving an expression for the iteration of operators
    $\bar P_\tau^n$.

    From the scheme \eqref{ex_scheme} and an iteration, we find that
    $$X_n=(1+b\tau)^nX_0+\sigma\sum_{k=0}^{n-1}(1+b\tau)^k \Delta B_{n-k}.$$
    Let  $\psi_n(\xi)$ and $\phi(\xi)=e^{-\frac12\xi^2}$ be the
    characteristic functions of $X_n$ and an $N(0,1)$ Gaussian variable, respectively.
    To  simplify the notations,
    let $C=C(b,\sigma,T,u_0)>0$ denote various positive constants
    depending only on $b$, $\sigma$, $T$, and $u_0$, and define
    \begin{align*}
      & \alpha:=  e^{bt}, \quad \alpha_n:=  (1+b\tau)^n; \quad
      \beta^2:= \frac{\sigma^2}{b}(e^{2bt}-1), \quad
      \beta_n^2:= \frac{\sigma^2}{b}\frac{2\big((1+b\tau)^{2n}-1\big)}{2+b\tau}.
    \end{align*}

    By the independence of the increments,
    the definition of $\Delta B_{k}$, and summation of exponents,
    \begin{align*}
        \psi_n(\xi)&=\psi_0(\alpha_n \xi)\prod_{k=0}^{n-1}
            \phi(\alpha_k\sqrt\tau\sigma\xi)
            = \psi_0(\alpha_n\xi)e^{- \frac 14 \beta_n^2\xi^2}.
    \end{align*}
    Note that $\psi_0(\xi)=\sqrt{2\pi}\hat u_0(\xi)$ and
    $\psi_n(\xi)=\sqrt{2\pi}{\hat{\bar u}}_n(\xi)$. Taking the
    inverse Fourier transform, and after elementary computations using
    convolution, scaling, and change of variables, we find that
    \begin{align*}
       \big(\bar P_\tau^n u_0\big)(y)& =\bar u_{n}(y)
       =\frac1{\sqrt{2\pi}}        \frac {u_0(\cdot/\alpha_n)}{ \alpha_n}
        * \left(\frac1{\sqrt{\beta_n^2/2}} e^{-\frac{\cdot^2}{ \beta_n^2}}\right)(y)\\
        &=\int_\R \frac{u_0(x)}{\sqrt{\pi\beta_n^2}}
        \exp\set{-\frac{(y- \alpha_nx)^2}{\beta_n^2}} dx.
    \end{align*}

    Now we estimate the $L^1$-norm of the difference of the exact and the
    approximate densities, $\|(\bar P_{\tau}^n-P_{t}) u_0\|_{L^1}$.
    Recall \eqref{conti_ker_lgv} and decompose
    \begin{align*}
       & \bar P_\tau^n u_0(y)-P_tu_0(y) \\
       & = \int_\R \frac{u_0(x)}{\sqrt{\pi\beta_n^2}}
        \exp\set{-\frac{(y- \alpha_nx)^2}{\beta_n^2}} dx
            -\int_\R \frac{u_0(x)}{\sqrt{\pi\beta^2}}
        \exp\set{-\frac{(y- \alpha x)^2}{\beta^2}} dx\\
       & = \beta_n\left(\frac 1{\beta_n}-\frac 1\beta\right)
            \int_\R \frac{u_0(x)}{\sqrt{\pi\beta_n^2}}
            \exp\set{-\frac{(y- \alpha_nx)^2}{\beta_n^2}} dx\\
       & \quad +\int_\R \frac{u_0(x)}{\sqrt{\pi\beta^2}}
            \left(\exp\set{-\frac{(y- \alpha_nx)^2}{\beta_n^2}}
                -\exp\set{-\frac{(y- \alpha_n x)^2}{\beta^2}}\right)dx\\
       & \quad + \int_\R \frac{u_0(x)}{\sqrt{\pi\beta^2}}
            \left(\exp\set{-\frac{(y- \alpha_n x)^2}{\beta^2}}
                -\exp\set{-\frac{(y- \alpha x)^2}{\beta^2}}\right)dx\\
       & =: E_1+E_2+E_3.
    \end{align*}
    We start with some observations.
    Let $a\in\R$. By l'H\^opital's rule and a Taylor expansion of $\ln(1+ax)$ about $x=0$,
    \[
        \lim_{x\to 0}\frac 1x \left( e^a-(1+ax)^{\frac 1x}\right)
            = \frac{a^2}{2}e^a,
    \]
    and hence there is a locally bounded function $\tilde C(a)>0$ such that
    \begin{equation}\label{eq est of approx e}
      \abs{e^a-\left( 1+ \frac an\right)^n}
        \leq \tilde C(a)\frac{a^2}{n}.
    \end{equation}

    Now for $E_1$, obviously $\norm{\int_\R \frac{u_0(x)}{\sqrt{\pi\beta_n^2}}
            \exp\set{-\frac{(\cdot- \alpha_nx)^2}{\beta_n^2}} dx}_1\leq \norm{u_0}_1$.
    We turn to analyse $\beta_n\left(\frac 1{\beta_n}-\frac 1\beta\right)$.
    Since
    \begin{equation}\label{eq property of e^z}
      \abs{e^x-1}\geq
      \begin{cases}
        \frac {\abs{x}}2,& -1\leq x,\\
        \frac 12,& x\leq -1,
      \end{cases}
    \end{equation}
    it follows that
    \[
      \beta^2:= \frac{\sigma^2}{b}(e^{2bt}-1)\geq
        \begin{cases}
            \sigma^2t,\qquad -1\leq 2bt,\\
            \frac{\sigma^2}{2|b|},\qquad 2bt\leq -1.
          \end{cases}
    \]
    In other words,
    \begin{equation}\label{eq est of 1/beta}
      \frac 1\beta\leq \frac{1\vee \sqrt{2|b|t}}{\sqrt{\sigma^2 t}}
        \quad\text{for all } t>0.
    \end{equation}
    Moreover
    \begin{equation}\label{eq est beta}
      \beta^2-\beta_n^2=\frac{\sigma^2}{b}\left(e^{2bt}-(1+bt/n)^{2n}\right)
            +\beta_n^2\frac{bt/n}{2+bt/n}.
    \end{equation}
    Hence combining \eqref{eq est beta} and \eqref{eq est of approx e}, we get
    \begin{align*}
      \abs{\beta_n\left(\frac 1{\beta_n}-\frac 1\beta\right)}
         =\abs{\frac{\beta^2-\beta_n^2}{\beta(\beta+\beta_n)}}
        \leq \abs{\frac{\beta^2-\beta_n^2}{\beta^2}}
        \leq \frac{\tilde C(2bt)4b^2t^2 }{n}\frac{\sigma^2}{2t}
        = \frac{Ct}{n}.
    \end{align*}
    Therefore, $\norm{E_1}_1\leq Ct/n$.

    Next, we estimate $E_2$. By a change of variable and
    the commuting property of convolutions, we have
    \begin{align*}
      E_2 & = \frac{1}{\alpha_n}\int_\R\frac{u_0\left( \frac{y-x}{\alpha_n} \right)}{\sqrt{\pi\beta^2}}
            \left( e^{-\frac{x^2}{\beta_n^2}}-e^{-\frac{x^2}{\beta^2}}\right)dx\\
      & = \frac{1}{\alpha_n}\int_\R\frac{u_0\left( \frac{y-\beta x}{\alpha_n} \right)}{\sqrt{\pi}}
            e^{-x^2}\left( e^{-\frac{\beta_n^2-\beta^2}{\beta_n^2}x^2}-1 \right)dx
    \end{align*}

    First we observe with \eqref{eq est of approx e}
    and \eqref{eq property of e^z} that
    \begin{align*}
      \beta_n^2 & \geq \frac{2\sigma^2}{3|b|} \abs{1-(1+b\tau)^{2n}}
        \geq \frac{2\sigma^2}{3|b|} \abs{1-e^{2bt}}- \abs{e^{2bt}-(1+b\tau)^{2n}}\\
       & \geq \frac {\sigma^2}3 \left(2t\wedge \frac 1{|b|}\right)-\frac{Ct^2}{n}
        \geq \frac {\sigma^2}6 \left(2t\wedge \frac 1{|b|}\right)\qquad
       \text{for}\qquad n>4CT(1\vee |b|T).
    \end{align*}
    Hence again from \eqref{eq est beta} and \eqref{eq est of approx e}  we have
    $\big|\frac{\beta_n^2-\beta^2}{\beta_n^2} x^2\big|\leq \frac{Ct}{n} x^2$.
    Therefore
    $$e^{-\frac{Ct}{n} x^2}-1\leq e^{\frac{\beta_n^2-\beta^2}{\beta_n^2}x^2}-1
        \leq e^{\frac{Ct}{n} x^2}-1.$$
    So  $\Big|e^{\frac{\beta_n^2-\beta^2}{\beta_n^2}x^2}-1\Big|
        \leq (1-e^{-\frac{Ct}{n} x^2})\vee (e^{\frac{Ct}{n} x^2}-1)
        = e^{\frac{Ct}{n} x^2}-1$.
   Now we can use Fubini's theorem to conclude that
    \begin{align*}
      \norm{E_2}_1 & \leq \frac{1}{\alpha_n}\iint
        \abs{\frac{u_0\left( \frac{y-\beta x}{\alpha_n} \right)}{\sqrt{\pi}}
            e^{-x^2}\left( e^{-\frac{\beta_n^2-\beta^2}{\beta_n^2}x^2}-1 \right)} dxdy \\
       & \leq   \norm{u_0}_1 \frac{1}{\sqrt{\pi}}
            \int_\R e^{-x^2}\big(e^{\frac{Ct}{n} x^2}-1\big)dx\\
       & = \norm{u_0}_1\left(\left(1-Ct/n\right)^{-1/2}-1\right)\\
       &  \leq \frac{Ct}{n}  \qquad \text{for } n \text{ large enough}.
    \end{align*}

    It remains to analyse $E_3$, and here the computations will differ
    according to the assumptions on $u_0$.
    We first assume that $\int_\R \abs{xu_0(x)}dx<\infty$.
    Denote $\gamma_n:=\frac{\alpha-\alpha_n}{\beta}$.
    Combining \eqref{eq est of approx e} and \eqref{eq est of 1/beta},
     it follows that
    \[
        \abs{\gamma_n}\leq \frac{\tilde C(bt)b^2t^2}{n}
            \frac{(1\vee \sqrt{2|b|t})}{\sqrt{\sigma^2 t}}
            \leq \frac{Ct^{3/2}}{n}.
    \]
    Therefore by a change of variable $y\mapsto y+\theta\gamma_n x$, we have
    \begin{align*}
      \norm{E_3}_1 & \leq \iint \abs{ \frac{u_0(x)}{\sqrt{\pi\beta^2}}
            \left(\exp\set{-\frac{(y- \alpha_n x)^2}{\beta^2}}
                -\exp\set{-\frac{(y- \alpha x)^2}{\beta^2}}\right)}dxdy \\
       & =  \iint \abs{\frac{u_0(x)}{\sqrt{\pi}}
            \left(e^{-(y+\gamma_n x)^2}
                -e^{-y^2}\right)}dydx\\
       & = \iint\frac{|u_0(x)|}{\sqrt{\pi}}\abs{\int_0^1 \gamma_n x\,
            2(y+\theta\gamma_n x)e^{-(y+\theta\gamma_n x)^2} d\theta  }dydx\\
       & \leq \frac{2|\gamma_n|}{\sqrt\pi}\int_\R \abs{xu_0(x)}dx
            \int_\R |z|e^{-z^2}dz\\
       & \leq \frac{Ct^{3/2}}{n}
             \qquad \text{for } n \text{ large enough}.
    \end{align*}

    For the other two cases,
    we write
    \begin{align*}
      E_3 & = \int_\R\left( \frac{1}{\alpha_n}u_0\Big(\frac{x}{\alpha_n}\Big)
            - \frac{1}{\alpha}u_0\left(\frac{x}{\alpha}\right)\right)
            \frac{1}{\sqrt{\pi\beta^2}} e^{-\frac{(y-x)^2}{\beta^2}}dx \\
      & = \left( \frac{1}{\alpha_n}- \frac{1}{\alpha}\right)
            \int_\R u_0\Big(\frac{x}{\alpha_n}\Big)
                 \frac{1}{\sqrt{\pi\beta^2}} e^{-\frac{(y-x)^2}{\beta^2}}dx\\
      & \quad + \frac 1\alpha \int_\R\left( u_0\Big(\frac{x}{\alpha_n}\Big)
                -u_0\Big(\frac{x}{\alpha}\Big)\right)
            \frac{1}{\sqrt{\pi\beta^2}} e^{-\frac{(y-x)^2}{\beta^2}}dx\\
      & =: E_3'+\frac 1\alpha E_3''.
    \end{align*}
    Note that $\alpha_n\to \alpha>0$ and hence $0<\alpha/2\leq\alpha_n\leq 2\alpha$ for large $n$.
    Since $\lim_n \alpha_n=\alpha$, and $\alpha,\alpha_n$ are bounded away
    from $0$. Combining  \eqref{eq est of approx e}, it follows that
    \[
      \abs{\frac 1\alpha-\frac 1{\alpha_n}}
        =\abs{\frac{e^{bt}-(1+bt/n)^n}{\alpha\alpha_n}}
        \leq \frac{1}{\alpha\alpha_n}\frac{\tilde C(bt)b^2t^2}{n}\leq \frac{Ct^2}{n}.
    \]
    Hence $\norm{E_3'}_1\leq \frac{Ct^2}{n}$ for large $n$.

    Regarding to $E_3''$, we assume now $\int_\R|xu_0'(x)|dx<\infty$. Then
    similarly
    \begin{align*}
      \norm{E_3''}_1  & \leq  \int_\R\abs{u_0\Big(\frac{x}{\alpha}\Big)
                -u_0\Big(\frac{x}{\alpha_n}\Big)}dx \\
       & \leq \int_\R \int_0^1 \abs{\left( \frac{1}{\alpha_n}- \frac{1}{\alpha}\right)x
            u_0'\left( \frac x\alpha +\theta \left( \frac{1}{\alpha_n}- \frac{1}{\alpha}\right)x\right)}d\theta dx\\
       & \leq C \abs{\frac{1}{\alpha_n}- \frac{1}{\alpha}}\int_\R|xu_0'(x)|dx
        \leq \frac{Ct^2}{n}.
    \end{align*}

    Finally we consider the most general case
    $u_0\in L^1(\R)$ only.
    We will prove a uniform convergence to zero without any rate.
    For any $\eps>0$, there exist
    $R>0$ and $N_1>0$ such that
    \[
        \int_{|x|>R}\abs{u_0\Big(\frac{x}{\alpha}\Big)}dx<\frac{ \eps}{3}\quad\text{and}\quad
            \int_{|x|>R}\abs{u_0\Big(\frac{x}{\alpha_n}\Big)}dx<\frac{ \eps}{3}\quad \text{for all } n>N_1.
    \]
    Next, since translation is continuous in $L^1$, i.e. $\lim_{h\to
      0}\int_{\R} |w(x)-w(x+h)|dx=0$ for any $w\in L^1(\R)$,
    and $h_n(x):=\left(\frac 1{\alpha_n}-\frac 1\alpha\right)x$ tends to $0$
    uniformly for $|x|\leq R$, there exists $N_2 >0$ such that for $n>N_2$,
    \begin{align*}
       \int_{|x|\leq R}\abs{u_0\Big(\frac{x}{\alpha}\Big)
                -u_0\Big(\frac{x}{\alpha_n}\Big)}dx
        = \int_{|x|\leq R}\abs{u_0\Big(\frac x\alpha\Big)
        -u_0\Big(\frac x\alpha+h_n(x)\Big)}dx
       < \frac{ \eps}{3}.
    \end{align*}
    Therefore we can conclude that $E_3''\to 0$ in $L^1$:
    \begin{align*}
      \norm{E_3''}_1 &
       \leq  \int_{|x|>R}\left(\abs{u_0\Big(\frac{x}{\alpha}\Big)}
            +\abs{u_0\Big(\frac{x}{\alpha_n}\Big)}\right)dx
            + \int_{|x|\leq R}\abs{u_0\Big(\frac{x}{\alpha}\Big)
                -u_0\Big(\frac{x}{\alpha_n}\Big)}dx \\
       & <   \eps \qquad \text{for}\qquad n>N_1\vee N_2.
    \end{align*}
    The proof is complete.


\bibliographystyle{amsplain}
\bibliography{Local-PI}

\end{document}